\documentclass[11pt]{article}
\usepackage{amsfonts,amsmath,amssymb,amsthm, amscd}
\usepackage[dvips]{epsfig,graphics}

\numberwithin{equation}{section}

\newtheorem{theorem}{Theorem}[section]
\newtheorem{prop}[theorem]{Proposition}

\newtheorem{cor}[theorem]{Corollary}

\theoremstyle{definition}
\newtheorem{definition}[theorem]{Definition}
\newtheorem{example}[theorem]{Example}
\newtheorem{remark}[theorem]{Remark}

\newcommand{\D}{\Delta}

\def\<{{\langle}}
\def\>{{\rangle}}

\def\Z{\mathbb Z}

\def\C{\mathbb C}

\def\ni{\noindent} 

\begin{document}

\title{Knot Group Epimorphisms, II}

\author{Daniel S. Silver \thanks{Partially supported by NSF grant
DMS-0304971.} \and Wilbur Whitten}

\maketitle 

\begin{abstract}   \noindent We consider the relations $\ge$ and $\ge_p$ on the collection of all knots, where $k \ge k'$ (respectively, $k \ge_p k'$) if there exists an epimorphism  $\pi k \to \pi k'$ of knot groups (respectively, preserving peripheral systems). When $k$ is a torus knot, the relations coincide and $k'$ must also be a torus knot; we determine the knots $k'$ that can occur. If $k$ is a $2$-bridge knot and $k \ge_p k'$, then $k'$ is a $2$-bridge knot with determinant a proper divisor of the determinant of $k$; only finitely many knots $k'$ are possible.  \end{abstract}

\begin{footnote}{Mathematics Subject Classification:  
Primary 57M25.}\end {footnote}

\section{Introduction} \label{Section 1} In recent years, numerous papers have investigated epimorphisms between knot groups and non-trivial maps between knot exteriors (or compact, orientable $3$-manifolds with boundary); see \cite{bb}, \cite{gr2}, \cite{ks}, \cite{ors}, \cite{kaw}, \cite{SWh1}, \cite{SWh2}, for example. We consider the first of these problems, which we formulate as follows (cf. \cite{ors}, \cite{riley1}). 

\medskip
{\sl 1. Given a nontrivial knot $k \subset {\mathbb S}^3$, 
classify the collections of knots $K$ for which there exists an epimorphism of knot groups $\pi K \to \pi k$, perhaps one preserving peripheral structure. 

2. For $k$ fixed, classify those knots $K$ for which there exists an epimorphism $\pi k \to \pi K$. }

\medskip

Let $k$ be a knot in ${\mathbb S}^3$, and let $E(k)$ denote the exterior of $k$. Orient both ${\mathbb S}^3$ and $k$. Choose and fix a point  $*$ on $\partial E(k)$, and set $\pi k = \pi_1({\mathbb S}^3\setminus k, *)$. Also, choose oriented curves $m$ and $l$ in $\partial E(k)$ meeting transversely at $*$ and representing a meridian-longitude system for $\pi k$; we use $m$ and $l$ to represent their classes in $\pi k$ as well. We consider knot group epimorphisms $\phi: \pi K \to \pi k$  defined up to automorphisms of $\pi k$. Since the type of a knot is determined by its complement \cite{gl}, an automorphism of $\pi k$ necessarily sends $m$ to a conjugate of $m$ or $m^{-1}$
\cite{tsau}. (We recall that knots have the {\it same type} if there exists an autohomeomorphism of ${\mathbb S}^3$ taking one knot to the other.) Hence we will call an element of $\pi k$ that is conjugate to $m$ or $m^{-1}$ a {\it meridian} of $k$, and we say that such an element is {\it meridional}.

Recall that if $k$ is nontrivial, then the inclusion-induced homomorphism $i_*: \pi_1(\partial E(k)) \to \pi_1(E(k))$ is injective
and defines a conjugacy class of subgroups of $\pi k$ -- the so-called
peripheral subgroups of $\pi k$, each member isomorphic to $\Z \times \Z$.  A homomorphism of knot groups {\it preserves peripheral structure} if it takes peripheral subgroups into peripheral subgroups. Recall also that for knots $k$ and $k'$ and epimorphism
$\phi: \pi k \to \pi k'$, we always have $\phi [(\pi k)'] = [(\pi k')']$, $\phi^{-1}[(\pi k')']=(\pi k)'$, and
${\rm ker} (\phi) \subset (\pi k)'$. Here $(\ )'$ denotes commutator subgroup. 

We write $k \ge k'$ whenever there exists an epimorphism $\phi: \pi k \to \pi k'$. If an epimorphism exists that preserves peripheral structure, then we write $k \ge_p k'$. The relation $\ge$ is a partial order on prime knots, while $\ge_p$ is a partial order on the collection of all knots \cite{SWh1}. (In \cite{SWh1} and \cite{SWh2} a slightly different notation is used.)  We write $k > k'$ if $k \ge k'$ but $k \ne k'$. The expression $k >_p k'$ has a similar meaning. 


In Section 2, we prove that if $\phi: \pi k \to \pi k'$ is an epimorphism taking a meridian of $k$ to a meridian of $k'$ and if $k'$ is prime, then $\phi$ preserves peripheral structure.  We prove several results about the relations $k \ge k'$ and $k\ge_p k'$ when $k$ is either a torus knot or a $2$-bridge knot. For a given torus knot $k$, Proposition \ref{Prop2} classifies those knots $k'$ for which there exists an epimorphism $\phi: \pi k \to \pi k'$, while Proposition \ref{Prop3} describes $\phi$ up to an automorphism of $\pi k'$. If $k$ is a $(p_1, q_1)$ $2$-bridge knot (with $p_1, q_1$ relatively prime odd integers, $p_1 \ge 3$ and $-p_1 < q_1 < p_1$) and if $k >_p k'$ with $k'$ nontrivial, then Proposition \ref{Prop4} asserts that $k'$ is a $(p_2, q_2)$ $2$-bridge knot such that $p_2$ properly divides $p_1$.

As a corollary to Proposition \ref{Prop4}, we show that given any $2$-bridge knot $k$, there are only finitely many knots $k'$ for which a meridian-preserving epimorphism $\pi k \to \pi k'$ exists. This is a partial answer to a problem of J. Simon (Problem 1.12 of \cite{kirby}). 
We close Section 2 with an example of two knots $k, k'$ for which there exists an epimorphism $\pi k \to \pi k'$ preserving meridians but taking the longitude of $k$ to $1$. Such epimorphisms correspond to zero-degree maps $E(k) \to E(k')$. 

In Section 3 we introduce the notions of minimal and $p$-minimal knots. We prove that twist knots are $p$-minimal, while a $(p_1, p_2)$-torus knot is minimal if and only if both $p_1$ and $p_2$ are prime.
Section 4 comprises a list of open questions. 

The second author thanks the Department of Mathematics at the University of Virginia for their continued hospitality and the use of their facilities.

\section{Epimorphisms and partial orders} \label{epimorphisms}

The followng proposition and its corollary give the useful fact that if $m$ is a meridian for a knot $k$, then an epimorphism $\phi: \pi k \to \pi k'$ such that $\phi(m)$ is a meridian of $k'$ also preserves peripheral structure, provided that $k'$ is a prime knot. 

\begin{prop} \label{Prop1} Let $k$ be a prime knot with  meridian-longtitude pair $(m, l)$ Then $Z(m) \cap (\pi k)''= \< l\>$
(= subgroup of $\pi k$ generated by $l$), where $Z(m)$ is the centralizer of $m$ in $\pi k$. \end{prop}

\begin{proof} Suppose that $g \in (\pi k)'', \ g\ne 1$, and $mg=gm$. If 
$g \in \< m, l\>$, then $g = l^d$, for some $d \ne 0$, since $g \in (\pi k)''$. We therefore assume that $g \notin \< m, l\>$. By Theorem 1 of \cite{simon}, $k$ is either a torus knot or a nontorus cable knot, since $k$ is prime. 

Assume first that $k$ is a torus knot, and set $P= \<m, l\>$. By Theorem 2 of \cite{simon}, $g^{-1}Pg \cap P$ is infinite cyclic (since
$mg=gm$). Since $g^{-1}Pg \cap P$ contains $m$, we have $g^{-1}Pg \cap P=\<m\>$. But $g^{-1}Pg \cap P$ also contains a generator of the center of $\pi k$. Since this is a contradiction, $k$ must be a nontorus, cable knot.

We have now that $E(k)= E(k_0)\cup_{T_0}S$, where $S$ is a cable space and $k$ is a cable about a nontrivial knot $k_0$. We can assume that $S$ is a component of the characteristic submanifold of $E(k)$. Note that $S$ is a small Seifert fibered manifold having an annulus with exactly one cone point as its base orbifold. Since $m$
and $l$ can be considered as elements of $\pi_1 S$ (well defined up to conjugation in $\pi k$), it follows from Theorem VI 1.6 (i) of \cite{js} that $Z(m)$ is a subgroup of $\pi_1 S$. Therefore, $g \in \pi_1 S$ (along with $m$ and $l$), and hence $g$ commutes with a generator of the center of $\pi_1 S$, which of course belongs to $P$.

As in the case of a torus knot, $g^{-1}Pg \cap P$ (as a subgroup of $\pi k$) is neither trivial nor infinite cyclic, which yields a contradiction. 
\end{proof}

\begin{remark}\label{rem2} 1. It is easy to see that the proposition does not hold if $k$ is composite. 

2. If $\phi: \pi k \to \pi k'$ is an epimorphism, then $\phi(l) \in Z(\phi(m))\cap (\pi k')''\ (= Z(\phi(m) \cap (\pi k')')$ \cite{jl}. In fact, given $k'$ and 
elements $\mu, \lambda \in \pi k'$, there exists a knot $k$ with meridian-longitude pair $(m, l)$ and an epimorphism $\phi: \pi k \to \pi k'$ such that $\phi(m)=\mu$ and $\phi(l)=\lambda$ if and only if $\mu$ normally generates $\pi k'$and $\lambda \in Z(\mu)\cap (\pi k')'' (see \cite{jl}).$ \end{remark}

\begin{cor}\label{Cor1.1} Let $k$ be a knot and $k'$ a prime knot. Let $(m,l)$ and 
$(m', l')$ be meridian-longitude pairs for $k$ and $k'$, respectively. If there exists an epimorphism $\phi: \pi k \to \pi k'$ with $\phi(m)= m'$, then $\phi$ preserves peripheral structure; in fact, $\phi(l) = (l')^d$, for some $d \in \Z$. \end{cor}

\begin{proof} As noted in Remark \ref{rem2} above, we must have $\phi(l) \in Z(\phi(m))\cap (\pi k')''\ (=Z(m')\cap (\pi k')'')$. Since $k'$ is prime, 
$Z(m') \cap (\pi k')'' = \<l\>$. Thus $\phi(l) = (l')^d$, for some $d\in \Z$.
\end{proof}

When we say that a knot is a $(p,q)$-torus knot, we will always assume that $p, q \ge 2$ and that $(p, q)=1$. Such a knot is necessarily nontrivial. 

\begin{prop}\label{Prop2} Let $k$ be a $(p_1, p_2)$-torus knot, and let $k'$ be a nontrivial knot. The following statements are equivalent. 

(1) $k \ge_p k',$

(2) $k \ge k',$

(3) $k'$ is an $(r_1, r_2)$-torus knot, for some $r_1, r_2 \ge 2$, such that $r_1|p_1$ and $r_2|p_2$, or $r_1|p_2$ and $r_2|p_1$.

\end{prop}

\begin{proof} Obviously, statement (1) implies statement (2).
Assume that $k \ge k'$. Then there exists an epimorphism $\phi: \pi k \to \pi k'$. If $k'$ is not a torus knot, then $\phi$ must kill the center of $\pi k$, since the only knots with groups having nontrivial centers are torus knots \cite{bz1}, and thus $\phi$ factors through the free product
$\Z_{p_1}*\Z_{p_2}$ of cyclic groups.  But no knot group can be a homomorphic image of $\Z_{p_1}*\Z_{p_2}$, since knot groups contain no nontrivial elements of finite order. Therefore, there exist integers $r_1, r_2\ge 2$ with $(r_1, r_2)=1$ such that $k'$ has the type of an $(r_1, r_2)$-torus knot. 

We have the following commutative diagram of epimorphisms
$$\begin{matrix} \pi k & \longrightarrow & \Z_{p_1}*\Z_{p_2} \\ \phi \downarrow & & \downarrow \psi \\  \pi k' & \longrightarrow & \Z_{r_1}*\Z_{r_2} \end{matrix}$$
in which the horizontal maps are canonical, and $\psi$ is the diagram-filling homomorphism. Let $t_1, t_2$ generate $\Z_{p_1}, \Z_{p_2},$
respectively. Since $\psi$ is an epimorphism, $\psi(t_1), \psi(t_2)$ generate $\Z_{p_1}*\Z_{p_2}.$ Moreover, each of $\psi(t_1)$ and $\psi(t_2)$ has finite order in $\Z_{p_1}*\Z_{p_2}.$ It follows from the torsion theorem for free products (see Theorem 1.6 of \cite{ls}, for example) that there are generators $s_1$ and $s_2$ of $\Z_{r_1}$ and $\Z_{r_2}$, respectively, such that either $\psi(t_1)= u_1 s_1 u_1^{-1}$ and 
$\psi(t_2) = u_2 s_2 u_s^{-1}$ or else $\psi(t_1) = u_1 s_2 u_1^{-1}$ and $\psi(t_2) = u_2 s_1 s_2^{-1}$,  for some $u_1, u_2 \in \Z_{p_1}*\Z_{p_2}$. Hence either $r_1|p_1$ and $r_2|p_2$ or else $r_1|p_2$ and $r_2|p_1$. 
Hence statement (2) implies statement (3). 

Finally, assume statement (3). Let $T$ be a standardly embedded, unknotted torus in ${\mathbb S}^3$ with complementary solid tori $V_1$ and $V_2$ such that $V_1 \cap V_2=T$. Assume that $C_i$ is an oriented core of $V_i$, for $i=1,2$, serving as an axis for periodic rotations of ${\mathbb S}^3$, each taking $T$ and the other axis to itself. Moreover, let $k$ be a $(p_1, p_2)$-torus knot in $T$ with $|{\rm lk}(k, C_i)|= p_i$, for $i=1, 2$, and such that periodic rotations of ${\mathbb S}^3$ of appropriate orders about each $C_i$ take $k$ to itself (see Proposition 14.27 \cite{bz2}). Assume that $r_1|p_1$ and $r_2|p_2$, and let $n_i r_i=p_i$, for $i=1,2$. A rotation of ${\mathbb S}^3$ of order $n_2$ about $C_1$ then yeilds a $(p_1, r_2)$-torus knot $k''$ as a factor knot. Similary, a rotation of ${\mathbb S}^3$ of order $n_1$ about the image axis of $C_2$ under the first rotation yields the $(r_1, r_2)$-torus knot $k'$ as a factor knot. Thus we have $k\ge_p k'' \ge_p k'$; that is, $k'$ is obtained from $k$ by at most two periodic rotations, each of which preserves peripheral structure. If $r_1|p_2$ and $r_2|p_1$, then the proof is similar.
\end{proof}

For a given torus knot $k$, Proposition \ref{Prop2} classifies those nontrivial knots $k'$ for which there exists an epimorphism $\phi: \pi k \to \pi k'$. The next result describes $\phi$ up to an automorphism of $\pi k'$. We recall from \cite{schreier} that an automorphism of the $(p,q)$-torus knot group $\<x, y \mid x^p = y^q\>$, with $p,q >1$ and $(p,q)=1$, has the form $x \mapsto w^{-1}x^\epsilon w,\ y\mapsto w^{-1}y^\epsilon w$, for $\epsilon \in \{-1, 1\}$.

\begin{prop}\label{Prop3} If $k$ and $k'$ are nontrivial torus knots with groups $\pi k =\< u, v \mid u^{p_1}= v^{p_2}\>$ and $\pi k' = \<a, b \mid a^{r_1} = b^{r_2}\>$ such that $r_i | p_i\  (i = 1, 2)$, and if $\phi: \pi k \to \pi k'$ is an epimorphism, then up to an automorphism of $\pi k'$, we have $\phi(u) = a^{n_2}$ and $\phi(v)= c^{-1} b^{n_1} c$, where $n_i r_i =p_i$ ($i=1, 2)$ and $c= b^s a^t$, for some $s, t \in \Z$.\end{prop}

\begin{proof} The element $(\phi(u))^{p_1}$ is in the center $Z(\pi k')$. Hence $\phi(u) = c_1^{-1}a^{\alpha_1}c_1$ or $\phi(u) = c_2^{-1}b^{\alpha_2}c_2$, for some $c_1, c_2 \in \pi k'$ and $\alpha_1, \alpha_2 \in \Z$ (see Lemma II. 4.2 \cite{js}). If $\phi(u)= c_2^{-1} b^{\alpha_2}c_2$ then $b^{\alpha_2 n_1 r_1}= b^{s r_2}$, for some $s$, and so $r_2|\alpha_2$, since $(r_2, n_1 r_1)=1$. But then $\phi(u) = c_2^{-1} b^{\alpha_2}c_2 \in Z(\pi k')$, which is a contradiction, since $\phi$ is an epimorphism. Thus $\phi(u) = c_1^{-1} a^{\alpha_1}c_1$ (and
$\phi(v) = c_2^{-1}b^{\alpha_2}c_2$, for some $c_2 \in \pi k'$ and $\alpha_2 \in \Z$).

Now $(\phi(u))^{p_1}= (\phi(v))^{p_2} \in \Z(\pi k')$, and so $a^{\alpha_1 p_1} = b^{\alpha_2 p_2}$; that is, $a^{r_1(n_1\alpha_1)}=b^{r_2(n_2 \alpha_2)}$. Since $a^{r_1}= b^{r_2}$ in $\pi k'$, it follows that $n_1 \alpha_1 = n_2 \alpha_2$; hence $n_d|\alpha_e$ for $d,e \in \{1,2\}$ and $d \ne e$, as $(n_1, n_2)=1$. Thus we can write
$\alpha_2 = n_1 \alpha_1 n_2^{-1}$ and get $\phi(u) = c_1 a^{\alpha_1}c_1$ and $\phi(v) = c_2^{-1} b^{n_1 \alpha_1 n_2^{-1}}c_2$,
where $\alpha_1$ is a multiple of $n_2$. Setting $\alpha_1 n_2^{-1} = n$, we have $\phi(u) = c_1^{-1} a^{n n_2}c_1$ and $\phi(v) = c_2^{-1} b^{n n_1} c_2$, for $0 < n \le \alpha_1.$

We show that $n=1$. Since $(p_1, p_2)=1$ and  $ip_1 + jp_2 = (in_1)r_1 + (j n_2)r_2 =1$, for some $i$ and $j$, the element $u^j v^i$ can be taken as a meridian of $k$ and the normal closure of $\phi(u^jv^i)$ is $\pi k'$. As a convenience, after conjugation of $\pi k'$ by $c_1$, we assume that $\phi(u) = a^{n n_2}$ and $\phi(v) = c^{-1} b^{n n_1} c$, where $c=c_2 c_1^{-1}$. So 
$$\phi(u^jv^i) = a^{(jn_2)n}c^{-1}b^{(in_1)n}c$$
$$\hskip 1.7 truein =(a^{(jn_2)n}b^{(in_1)n})(b^{-(in_1)n}c^{-1}b^{(in_1)n}c).$$
Since $\phi(u^j v^i)$ normally generates $\pi k'$, we have $|{\rm lk}(k',m )|=1$, where $m$ represents $\phi(u^j v^i)$. Since $\phi(u^jv^i)$
and $a^{(jn_2)n}b^{(in_1)n}$ have the same abelianizations, this linking number is $r_1(in_1)+r_2(j n_2)$. Thus $\phi(u)= a^{n_2}$ and $\phi(v) = c^{-1} b^{n_1}c$.

Now $A= \{a^{n_2}, c^{-1}b^{n_1}c\}$ generates $\pi k'$ (by assumption), and since $B = \{a, c^{-1}b^{n_1}c\}$ generates $A$, then $B$ generates $\pi k'$. Similarly, $C=\{a, c^{-1}b c\}$ generates $B$ and hence $\pi k'$. Taking $C$ as the generating set of $\pi k'$, it is now an exercise to show that $c$ (in $\Z_{r_1}*\Z_{r_2}$) has the form $b^s a^t$ ($1 \le s \le r_2, 1 \le t \le r_1$).  (In fact, such an exercise appears as Excercise 15, page 194, of \cite{mks}.) Thus $\psi^{-1}(b^sa^t) = (b^sa^t)Z(\pi k')$, where $\psi: \< a, b \mid a^{r_1} = b^{r_2}\> \to \< a, b \mid a^{r_1}, b^{r_2}\> = \Z_{r_1}*\Z_{r_2}$ is defined by $a\mapsto a$ and $b \mapsto b$ so that $\{ a, c^{-1}bc\}$ generates $\Z_{r_1}*\Z_{r_2}$.

\end{proof}

\begin{cor}\label{Cor3.1} Torus-knot group epimorphisms preserve peripheral structure. \end{cor}

\begin{proof} Let $k$ be a $(p_1, p_2)$-torus knot, and let $k'$ be a nontrivial knot. Suppose that there exists an epimorphism $\phi: \pi k \to \pi k'$. By Proposition \ref{Prop2}, $k'$ is an $(r_1, r_2)$-torus knot, and we can assume that $n_i r_i = p_i$, $(i=1,2)$. We have 
$\pi k =\< u,v \mid u^{p_1} = v^{p_2}\>, \ \pi k'= \< a, b \mid a^{r_1} = b^{r_2}\>$, and $ip_1+ jp_2 = (i n_1)r_1 + (j n_2) r_2 = 1,$ for some $i, j$. The element $u^j v^i$ is a meridian of $k$, and according to Proposition \ref{Prop3}, we can assume that $\phi(u) = a^{n_2}$ and $\phi(v) = c^{-1} b^{n_1} c$, where $c = b^s a^t$ ($s, t \in \Z)$. Thus
$u^j v^i \mapsto a^{j n_2}(a^{-t}b^{-s})b^{in_1}(b^sa^t)= a^{-t}(a^{jn_2}b^{in_1})a^t$, which is clearly a meridian. It follows from Corollary \ref{Cor1.1} that $\phi$ preserves peripheral structure. In fact, if $l_1$ and $l_2$ are the (appropriate) longitudes of $k$ and $k'$, respectively, then $\phi(l_1) = a^{-t} l_2^{n_1 n_2} a^t$. \end{proof}

\begin{cor} \label{Cor3.2} If $k$ is a torus knot and if $k\ge k'$, then $\pi k'$ embeds in $\pi k$. \end{cor}

\begin{proof} If $k$ is a $(p_1, p_2)$-torus knot, then $k'$ is an $(r_1, r_2)$-torus knot, for some $r_1, r_2 \ge 2$ with $(r_1, r_2)=1$, and either $r_1|p_1$ and $r_2|p_2$ or else $r_1|p_2$ and $r_2|p_1$. It follows immediately that $\pi k'$ embeds in $ \pi k$ (see Theorem 5.1 \cite{gw}). 
\end{proof}

\begin{remark} 1. If $k$ is a $(p_1, p_2)$-torus knot, then there may well exist an $(r_1, r_2)$-torus knot $k'$ such that $\pi k'$ embeds in $\pi k$ but it is not the case that
$k \ge k'$. For example, let $p_1 =2$ and $p_2 = 3 \cdot 5$, and take $r_1 =3, r_2= 5$. Then $\pi k'$ embeds in $\pi k$ by \cite{gw}, but 
it is not the case that $k \ge k'$ by Proposition \ref{Prop2}. 

2. By Corollary \ref{Cor3.2}, we know that if $k$ is a torus knot, then $k\ge k'$ implies that $\pi k'$ is a subgroup of $\pi k$. The index of this embedding is finite, however, if and only if $k' = k$ (see Remark 3, page 42 of \cite{gw}). \end{remark}

\begin{cor} \label{Cor3.3} If $k$ is a torus knot and $k \ge k'$, then the crossing number of $k$ is no less than that of $k'$. \end{cor}

\begin{proof} This follows from Proposition \ref{Prop2} and the fact
that   \cite{mur} the crossing number of a $(p,q)$-torus knot is $\min \{ p(q-1), q(p-1)\}$. \end{proof} 

We are particularly interested in the relative strengths of the two relations $\ge$ and $\ge_p$. Stated in general terms, our inquiry takes the form: \bigskip

Q1. Given knots $k$ and $k'$, when does $k \ge k'$ imply $k \ge_p k'$? \bigskip

\ni For knots $k$ and $k'$ with at most 10 crossings $k \ge k'$ implies 
$k \ge_p k'$ by \cite{ks}. Question 1 generates a number of related questions. One of them is: \bigskip

Q2. For which pairs of knots $k$ and $k'$ does there exist an epimorphism $\pi k \to \pi k'$ but no epimorphism preserving meridians? \bigskip

For the present, we will consider the case $k \ge_p k'$ with $k$ a 2-bridge knot. When we say that $k$ is a $(p, q)$ 2-bridge knot, we assume that $p, q$ are relatively prime odd integers, $p \ge 3$ and $-p < q < p$. Recall that $p$ is $\det (k)$, the determinant of $k$. 

A representation $\pi \to {\rm SL}_2{\mathbb C}$ is {\it parabolic} if if projects to a parabolic representation $\pi \to {\rm PSL}_2{\mathbb C}={\rm SL}_2{\mathbb C}/\<-I\>$ sending some and thus every meridian to a parabolic element.

\begin{prop} \label{Prop4} Let $k$ be a $(p_1, q_1)$ 2-bridge knot and let $k'$ be a nontrivial knot. If $k >_p k'$, then $k'$ is a $(p_2, q_2)$ 2-bridge knot such that $p_2$ properly divides $p_1$ (and hence $\D_{k_2}(t)$ properly divides $\D_{k_1}(t)$.) \end{prop}

\begin{proof} Let $m_1$ be a meridian of $k$ and $m_2$ a meridian of $k'$. Since $k \ge_p k'$, we have an epimorphism $\phi: \pi k \to \pi k'$ with $\phi(m_1) = m_2$, which induces an epimorphism
$\pi k/\<\< m_1^2\>\> \to \pi k'/\<\< m_2^2\>\>$ of $\pi$-orbifold groups. Since $k$ is a $(p_1, q_1)$ 2-bridge knot, $\pi k/\<\<m_1^2\>\>$ is the dihederal group $D_{p_1}$ of order $2p_1$. Hence $\pi k'/\<\<m_2^2\>\>$ is (finite) dihedral or $\Z_2$. By the Smith Conjecture 
\cite{bm}, the group $\pi k'/\<\< m_2^2\>\>$ is not $\Z_2$ (since by hypothesis,
$k'$ is not trivial) and thus it is isomorphic to $D_{p_2}$, for some $p_2$. Hence $p_2$ divides $p_1$, and therefore $p_2$ is odd. It follows that $k'$ is a $(p_2, q_2)$ 2-bridge knot, for some $q_2$; see Proposition 3.2 of \cite{bz}. Note that Proposition 3.2 of \cite{bz} depends on Thurston's orbifold geometrization theorem; see \cite{blp}, for example).

To see that $p_1 > p_2$, we examine two cases. First assume that $k$ is a 2-bridge torus knot (a $(p_1, 2)$-torus knot). By Proposition 
\ref{Prop2}, $k'$ is a $(p_2, 2)$-torus knot and $p_1 > p_2$, since $k >_p k'.$

For the second case, we assume that $k$ is hyperbolic, and we apply Riley's parabolic representation theory \cite{riley1}. Accordingly, if $k$ is a $(p, q)$ 2-bridge knot, then there are exactly $(p-1)/2$ conjugacy classes of nonabelian parabolic ${\rm SL}_2{\mathbb C}$ representations, corresponding to the roots of a monic polynomial $\Phi_{p,q}(w)$. As $\phi: \pi k \to \pi k'$ preserves peripheral structure, each parabolic representation $\theta': \pi k' \to {\rm SL}_2{\mathbb C}$ induces a parabolic representation  $\theta: 
\pi k  \to {\rm SL}_2{\mathbb C}$, and since $\phi$ is an epimorphism, $\phi$ induces a one-to-one function of conjugacy classes. When $p_2 = p_1$, the function is a bijection, and hence some representation $\theta'$ must induce an injection $\theta: \pi k \to {\rm SL}_2{\mathbb C}$, a lift of the faithful discrete representation $\pi k \to {\rm PSL}_2{\mathbb C}$ corresponding to the hyperbolic structure of ${\mathbb S}^3\setminus k$ (see \cite{thurston}). Since $\theta = \theta' \circ \phi$, the epimorphism $\phi$ is in fact an isomorphism, a contradiction as $k$ and $k'$ have different types. Hence $p_1 > p_2$. 

From the fact that $p_2 = |\D_{k'}(-1)|$ and $p_1 = |\D_k(-1)|$, it follows that $\D_{k'}(t)$ properly divides $\D_k(t)$. 
\end{proof}

The following corollary provides a partial answer to a problem of J. Simon (see Problem 1.12 of \cite{kirby}). 

\begin{cor}\label{simon} Let $k$ be a 2-bridge knot. There exist only finitely many knots $k'$ for which a meridian-preserving epimorphism $ \pi k \to \pi k'$ exists. \end{cor}

\begin{proof} Assume that a meridian-preserving epimorphism $\phi: \pi k \to \pi k'$ exists. Since $\pi k$ is generated by two elements, the same is true of $\pi k'$. By \cite{norwood}, $k'$ is 
a prime knot. Corollary \ref{Cor1.1} implies that $\phi$ preserves peripheral systems. By Proposition \ref{Prop4}, the knot $k'$ is 2-bridge.
The Alexander polynomial of $k'$ must divide that of $k$, and by 
\cite{riley2}, only finitely many possible such knots $k'$ exist. \end{proof}

\begin{remark}\label{rileypolys} Given a $(p, q)$  2-bridge knot $k$, one can use the Riley polynomial $\Phi_{p,q}$ to determine all knots $k'$ such that $k \ge_p k'$. Properties and applications of Riley polynomials will be discussed in a forthcoming paper. \end{remark}

\begin{cor} \label{nonvanishinglongitude} If $k$ is a nontrivial $2$-bridge knot, then any meridian-preserving epimorphism $\phi: \pi k \to \pi k'$ maps the longitude nontrivially. \end{cor}

\begin{proof}  By Proposition \ref{Prop4}, $k'$ is a $2$-bridge knot. 
Assume that $\phi$ maps the longitude of $k$ trivially. Let $\theta': \pi k' \to {\rm SL}_2 \C$ be any nonabelian parabolic representation. Then $\theta' \circ \phi$ is a nonabelian parabolic representation of $\pi k$ sending the longitude to the identity matrix, contradicting Lemma 1 of \cite{riley3}. Thus $\phi$ maps the longitude of $k$ nontrivially. 

\end{proof}
\begin{remark} \label{3remarks}

1. In \cite{ors}, the authors give a sufficient condition for the existence of a peripheral-structure preserving epimorphism between 2-bridge link groups. The condition is in fact a very efficient machine for generating many such epimorphisms. For example, one can use it to show that $k >_p k'$ for $k$ the $(175, 81)$ 2-bridge knot and $k'$ the $(7,3)$ 2-bridge knot, as pointed out by K. Murasugi. 

2. F. Gonzalez-Acu\~na and A. Ramirez \cite{gr2} proved that $k \ge_p \tau_{a,b}$, where $\tau_{a, b}$ is some torus knot, if and only if $k$ has property Q. (A knot $k$ has property Q if there is a closed surface $F$ in ${\mathbb S}^3 = X  \cup_F Y$ such that $k \subset F$ and $k$ is imprimative in each of $H_1(X)$ and $H_1(Y)$. Basic examples of such knots are torus knots.) In \cite{gr1}, they determined the 2-bridge knots $k$ such that $k \ge_p \tau_{a, 2}$ for some odd $a \ge 3$. 

3. Define a {\it knot manifold} to be a compact, connected, orientable, irreducible 3-manifold with boundary an incompressible torus. Such a manifold is said to be {\it small} if it contains no closed essential surface. A 3-manifold $M$ {\it dominates} another 3-manifold $N$ if there is a continuous, proper map $f: M \to N$ of nonzero degree.
(Here proper means that $f^{-1}(\partial N) = \partial M$.) A knot manifold is {\it minimal} if it dominates only itself. In \cite{bb}, Boileau and Boyer show that twist knots and $(-2, 3, n)$-pretzel knots ($n$ not divisible by 3) are minimal.  \end{remark} 

Suppose that $k >_p k'$ with $k'$ nontrivial, and let $\phi: \pi k \to \pi k'$ be an epimorphism that preserves peripheral structure. If $(m_1, l_1)$ and $(m_2, l_2)$ are fixed meridian-longitude pairs for $\pi k$ and $\pi k'$, respectively, then we can assume that $\phi(m_1) = m_2$ and $\phi(l_1) = l_2^d$, for some $d\in \Z$. Then $\phi$ is induced by a proper map $f: E(k) \to E(k')$, and the absolute value of the degree of $f$ is $|d|$, since $f_*: H_3(E(k), \partial E(k)) \to H_3(E(k'), \partial E(k'))$ takes the top class of $H_3(E(k), \partial E(k))$ to ${\rm deg} (f)$ times the top class of $H_3(E(k'), \partial E(k'))$ (see Proposition 6.2 of \cite{ors}). 



\begin{example} From the proof of Corollary \ref{Cor3.1}, it follows that an epimorphism $\pi k \to \pi k'$ in which each of $k$ and $k'$ is a nontrivial torus knot is induced by a nonzero-degree map; that is, a map sending $m_1 \mapsto m_2$ and $l_1 \mapsto l_2^d$ for some nonzero $d$. It is easy, however, to find knots $k$ and $k'$ and an epimorphism $\pi k \to \pi k'$ with $m_1 \mapsto m_2$ but $l_1 \mapsto 1$. For this, one can choose $k$ to be the square knot and $k'$ the trefoil. 


As a second example, let $R_y$ denote the knot of \cite{riley1} termed his ``favorite knot." A surgery description of $R_y$ appears in Figure 1.
It is not difficult to find an epimorphism $\pi R_y \to \pi (k_2 \sharp k_3)/\<\<zb^{-1}\>\>,$ where $k_2 \sharp k_3$ is the square knot indicated in Figure 1. Since the longitude of $k_2 \sharp k_3$  goes to $1$ (as does $zb^{-1}$) under the appropriate epimorphism $\pi (k_2 \sharp k_3) \to \pi 3_1$, we have a meridian-preserving, longitude-killing epimorphism $\pi R_y \to \pi3_1$. A similar argument shows that $8_{20} > 3_1$ (with longitude sent to $1$), but neither is $R_y > 8_{20}$ nor is $8_{20} > R_y$, since $R_y$ and $8_{20}$ are prime fibered knots of the same genus. Other such examples can be found in \cite{gr2}. 

\begin{figure}
\begin{center}
\includegraphics[height=2.5 in]{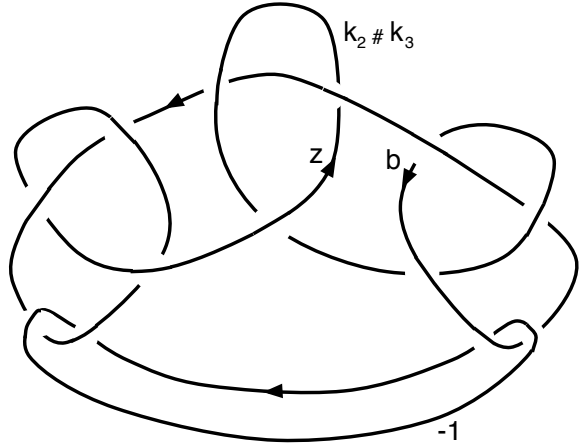}
\caption{Surgery description of Riley's knot}
\label{Figure1}
\end{center}
\end{figure}



\end{example}

\section{Minimality} \label{minimality}

\begin{definition} 1. A knot $k$ is {\it minimal} if $k \ge k'$ implies that $\pi k \cong \pi k'$ or else $k'$ is trivial. 

2. $k$ is $p$-{\it minimal} if $k \ge_p k'$ implies that $k=k'$ or else $k'$ is trivial. \end{definition}

\ni (Compare these definitions of ``minimality" with that given in the second part of Remark \ref{3remarks} above.)

Recall that the Alexander polynomial of a nontrivial 2-bridge knot is not equal to $1$. The next result follows immediately from Proposition \ref{Prop4}. 

\begin{cor} \label{Cor4.1} If $k$ is a nontrivial $(p_1, q_1)$ 2-bridge knot and if no proper nontrivial factor of $\D_k(t)$ is a knot polynomial, then $k$ is p-minimal. In particular, $k$ is p-minimal if $\D_k(t)$ is irreducible or if $p_1$ is prime. \end{cor}

\begin{remark} A nontrivial $(p_1, q_1)$ 2-bridge knot $k$ can have both $\D_k(t)$ reducible and $p_1$ composite but still be p-minimal. The simplest example is $k= 6_1$, the $(9,4)$ 2-bridge knot. \end{remark}

\begin{cor} \label{Cor4.2} 1. Every nontrivial  twist knot is p-minimal.

2. For each $n\ge 3$, there is a p-minimal knot with crossing number $n$. 
\end{cor}

\begin{proof} In view of Corollary \ref{Cor4.1}, we need only note that the Alexander polynomial of a nontrivial twist knot is quadratic and that for each $n\ge 3$ there is a twist knot with crossing number $n$. \end{proof}

\begin{cor} \label{Cor4.3} Every genus-one 2-bridge knot is p-minimal. \end{cor}

\begin{proof} Let $k$ be a genus-one 2-bridge knot. Since $k$ is alternating, the degree of its Alexander polynomial $\D_k(t)$ is 2. 
If $k'$ is a nontrivial knot and $k >_p k'$, then $k'$ is a 2-bridge knot
by Proposition \ref{Prop4}, and $\D_{k'}(t)$ divides $\D_k(t)$ properly. Since the degree of a knot polynomial is even, the degree of $\D_{k'}(t)$ must be 0. This is impossible, however, since $k'$ is nontrivial and 
alternating. 

\end{proof}

From Proposition \ref{Prop2} we have:

\begin{cor} \label{Cor4.4} If $k$ is a $(p_1, p_2)$-torus knot, then $k$ is minimal if and only if both $p_1$ and $p_2$ are prime. \end{cor}

Corollaries  \ref{Cor4.1}-\ref{Cor4.4} should be compared with theorems 3.16, 3.19 and corollaries 3.17, 3.18, 3.20 of \cite{bb}

\section{Questions.} If $k \ge_p k'$, then what properties of $k'$ can we deduce from $k$? For example,  if $k$ is fibered, then so is $k'$.

Not all properties of $k$ are inherited by $k'$. For example, if $k$ is prime, then $k'$ need not be \cite{SWh1}.  

Q3: If $k$ is alternating, must $k'$ also be alternating? 
[Yes, if $k$ is 2-bridge. This follows from Proposition \ref{Prop4} together with the fact that $2$-bridge knots are alternating \cite{bs}.]

Q4: Must the genus of $k'$ be less than or equal to that of $k$? [Yes, if $k$ is a 2-bridge knot or a fibered knot. In these cases $k'$ is 
a knot with the same property and  $\D_{k'}(t)$ divides $\D_k(t)$.
However, the genus of a $2$-bridge or fibered knot  equal to half the degree of its Alexander polynomial \cite{crowell}, \cite{bz2}. See also Proposition 3.7 of \cite{SWh1}.]

Q5: Must the crossing number of $k'$ be less than or equal to that of $k$? [Yes, if $k$ is a torus knot. This follows from Corollary \ref{Cor3.3}.]

Q6: Must the Gromov invariant of $k'$ be less than or equal to that of $k$? [Yes, if $k$ is a 2-bridge knot or a torus knot. If $k$ is a $2$-bridge knot, then any nontrivial epimorphism $\phi: \pi k \to \pi k'$ maps the longitude of $k$ nontrivially, by Corollary \ref{nonvanishinglongitude}. There exists a map $E(k) \to E(k')$ of degree $d>0$, and hence the Gromov invariant of $k$ is at least $d$ times that of $k'$. If $k$ is a torus knot, then so is $k'$, by Proposition \ref{Prop2}. Both $k$ and $k'$ have vanishing Gromov invariant \cite{soma}.] 


 \bigskip

\noindent {\sl Addresses:} Department of Mathematics and  Statistics, ILB 325, University of South Alabama, Mobile AL  36688 USA \medskip
\noindent 1620 Cottontown Road, Forest VA 24551 USA

\noindent {\sl E-mail:} silver@jaguar1.usouthal.edu; BJWCW@aol.com

\begin{thebibliography}{1}

\bibitem{bs} C. Bankwitz and H.G. Schumann, {\it \"Uber Viergeflechte}, {\sl Abh. Math. Sem. Univ. Hamburg\ \bf10} (1934), 263--284.

\bibitem{bb} M. Boileau and S. Boyer, {\it On character varieties, sets of discrete characters, and non-zero degree maps}, arXiv: GT/070138 v1.

\bibitem{blp} M. Boileau, B. Leeb and J. Porti, {\it Geometrization of 3-dimensional orbifolds}, {\sl Ann. Math.\ \bf162} (2005), 195--290.


\bibitem{bm} H. Bass and J. Morgan, The Smith Conjecture, Academic Press, New York, 1984.

\bibitem{bz} M. Boileau and B. Zimmermann, {\it The $\pi$-orbifold group of a link}, {\sl Math. Z.\ \bf200}, (1989), 187--208.


\bibitem{bz1} G. Burde and H. Zieschang, {\it Eine Kennzeichnung der Torusknoten}, {\sl Math. Ann. \ \bf 167} (1966), 169--176.

\bibitem{bz2} G. Burde and H. Zieschang, Knots, De Gruyter, Berlin, 1985.

\bibitem{crowell} R.H. Crowell, {\it Genus of alternating link types}, {\sl Ann. of Math.\ \bf3} (1969), 258--275.


\bibitem{gr2} F. Gonzalez-Acu\~na and A. Ramirez, {\it Epimorphisms of knot groups onto free products}, {\sl Topology \bf42} (2003), 1205--1227. 

\bibitem{gr1} F. Gonzalez-Acu\~na and A. Ramirez, {\it Two-bridge knots with property Q}, {\sl Quart. J. Math.\ \bf52}  (2001), 447--454.

\bibitem{gw} F. Gonzalez-Acu\~na and W. Whitten, Imbeddings of three-manifold groups, Memoirs Amer. Math. Soc. {\bf 99} (1992), No. 474. 


\bibitem{gl} C. Gordon and J. Luecke, {\it Knots are determined by 
their complements}, {\sl J. Amer. Math. Soc.\ \bf 2} (1989), 371--415.

\bibitem{js} W. Jaco and P. Shalen, {\it Seifert fibered spaces in $3$-manifolds}, {\sl Memoirs of the Amer. Math. Soc.\ \bf 21} (1979).


\bibitem{jl} D. Johnson and C. Livingston, {\it Peripherally specified homomorphs of knot groups}, {\sl Trans. Amer. Math. Soc.\ \bf 311} (1989), 135--146.

\bibitem{ls} R.C. Lyndon and P.E. Schupp, Combinatorial Group Theory, Springer-Verlag, Berlin, 1977. 

\bibitem{kaw} A. Kawauchi, {\it Almost identical imitations of $(3,1)$-dimensional manifolds}, {\sl Osaka J. Math. \bf 26} (1989), 743--758.

\bibitem{kirby} R. Kirby, Problems in low-dimensional topology, Geometric topology, Edited by H.Kazez,
AMS/IP, Vol. 2, International Press, 1997.

\bibitem{ks} T. Kitano and M. Suzuki, {\it A partial order on the knot table}, {\sl Experimental Math. \bf 14} (2005), 385--390. 

\bibitem{mks} W. Magnus, A. Karrass and D. Solitar,  Combinatorial Group Theory,  Wiley, New York, 1966. 

\bibitem{mur} K. Murasugi, {\it On the braid index of alternating links}, {\sl Trans. Amer. Math. Soc.\ \bf123} (1991), 237--260. 

\bibitem{norwood} F. H. Norwood, {\it Every two-generator knot is prime}, {\sl Proc. Amer. Math. Soc.\ \bf86} (1982), 143--147.

\bibitem{ors} T. Ohtsuki, R. Riley and M. Sakuma, {\it Epimorphisms between 2-bridge link groups}, preprint.

\bibitem{riley1} R. Riley, {\it Parabolic representations of knot groups, I}, {\sl Proc. London Math. Soc.\ \bf24} (1972), 217--242.

\bibitem{riley2} R. Riley, {\it A finiteness theorem for alternating links},
{\sl J. London Math. Soc. \bf 5} (1972), 263Ð266.

\bibitem{riley3} R. Riley, {\it Knots with parabolic property P}, {\sl Quart. J. Math.\ \bf25} (1974), 273--283.



\bibitem{schreier} O. Schreier, {\it Uber die gruppen $A^aB^b=1$},
{\sl Abh. Math. Sem. Univ. Hamburg\ \bf 3} (1924), 167--169.

\bibitem{SWh1} D.S. Silver and W. Whitten, {\it Knot group epimorphisms}, {\sl J. Knot Theory Ramifications\ \bf 15} (2006), 153--166. 


\bibitem{SWh2} D.S. Silver and W. Whitten, {\it Hyperbolic covering knots}, {\sl Alg. Geom. Topol.\ \bf 5} (2005), 1451--1469. 

\bibitem{simon} J. Simon, {\it Roots and centralizers of peripheral elements of knot groups}, {\sl Math. Ann.\ \bf 229} (1976), 205--209.

\bibitem{SnapPea} J. Weeks, www.geom.umin.edu

\bibitem{soma} T. Soma, {\it The Gromov invariant of links}, {\sl Invent. Math.\ \bf64} (1981), 445--454.

\bibitem{thurston} W.P. Thurston, The Geometry and Topology of 3-Manifolds, Princeton University, 1980. 

\bibitem{tsau} C.M. Tsau, {\it Algebraic meridians of knot groups}, {\sl Trans. Amer. Math. Soc.\ \bf 294} (1986), 733--747.





\end{thebibliography}
\end{document}